\documentclass[]{elsarticle}
\usepackage{setspace}
\doublespace
\usepackage{amsmath}
\usepackage{amssymb}
\usepackage{amsfonts}
\usepackage{subfigure}
\usepackage{graphicx}
\usepackage[british]{babel}
\usepackage{german}
\usepackage{url}
\newcommand{\ch}{{\binom{X}{2}}}

\usepackage{setspace}
\newcommand {\junk}[1]{}%
\usepackage{epic}
\usepackage{latexsym}
\usepackage{amsmath}
\usepackage{amssymb}

\usepackage[mathscr]{euscript}

\usepackage[all]{xy}

\usepackage{yfonts}

\DeclareMathAlphabet{\mathpzc}{OT1}{pzc}{m}{it}

\newtheorem{lemma}{Lemma}[section]
\newtheorem{theorem}{Theorem}
\newtheorem{proposition}[lemma]{Proposition}
\newtheorem{corollary}[lemma]{Corollary}

\newcommand{\0}{{\emptyset}}
\newcommand{\rk}{{\rm rk}}
\newcommand{\res}{{\rm res}}
\newcommand{\Hom}{{\rm Hom}}
\newcommand{\rr}{{\mathbb R}}

\newcommand{\cp}{{\mathcal P}}

\newcommand{\rG}{{\mathbb G}}
\newcommand{\rI}{{\mathbb I}}
\newcommand{\rC}{{\mathbb C}}
\newcommand{\rM}{{\mathbb M}}
\newcommand{\rB}{{\mathbb B}}
\newcommand{\sgn}{{\rm sgn}}
\newcommand{\supp}{{\rm supp}}

\newcommand{\cl}{{\mathcal L}}

\newcommand{\nn}{{\mathbb N}}

\newcommand{\RR}{{\mathbb R}}

\newcommand{\ra}{{\rightarrow}}
\newcommand{\LRa}{{\Leftrightarrow}}
\newcommand{\Ra}{{\Rightarrow}}
\newcommand{\pf}{\noindent{\em Proof: }}
\newcommand{\epf}{\hfill\hbox{\rule{3pt}{6pt}}\\}
\begin{document}

\title{A matroid associated with a phylogenetic tree}
\author{Andreas W.M.~Dress, Katharina T. Huber, and Mike Steel}
\address{(A.W.D.): CAS-MPG Partner Institute and Key Lab for Computational Biology, 320 Yue Yang Road, 200031 Shanghai, China;\\
(K.T.H.): School of Computing Sciences,
University of East Anglia, Norwich, UK;\\
(M.S.): Department of Mathematics and Statistics,
University of Canterbury, Christchurch, New Zealand.
}


\begin{abstract}
A (pseudo-)metric $D$ on a finite set $X$ is said to be a `tree metric' if there is a finite tree with leaf set $X$ and non-negative edge weights so that, for all $x,y \in X$, $D(x,y)$ is the path distance in the tree between $x$ and $y$.  It is well known that not every metric is a tree metric.
However, when some such tree exists, one can always find one whose interior edges have strictly positive edge weights and that has no vertices of degree $2$, any such tree is  -- up to canonical isomorphism -- uniquely determined by $D$, and one does not even need all of the distances in order to fully (re-)construct the tree's edge weights in this case.
Thus, it seems of some interest to investigate 
which subsets of $\binom{X}{2}$ suffice
to determine (`lasso')  
these edge weights.
In this paper, 
we use the results of a previous paper 
 to discuss the structure of 
a matroid that can be associated with an (unweighted) $X-$tree $T$ 
defined by the requirement that its
bases are exactly the `tight edge-weight lassos' for $T$, i.e, the minimal subsets $\cl$ of $\ch$ that lasso the edge weights of $T$.
\end{abstract}

\maketitle

\newpage

\section{Introduction}\label{intro}

Given any 
finite tree $T$ without vertices of degree $2$, there is an associated matroid $\rM(T)$ having ground  set $\binom{X}{2}$
 where $X$ is the set of leaves of $T$.  In this paper, we describe this matroid and investigate a number of interesting properties it exhibits. The motivation for studying this matroid is its relevance to the problem of uniquely reconstructing an edge-weighted tree from 
 its topology and just some
of the leaf-to-leaf distances in  that
tree. This combinatorial problem arises in phylogenetics (the inference of evolutionary relationships from genetic data) since -- due to patchy taxon coverage by available genetic loci \cite{patchy} -- reliable estimates of evolutionary distances 
can often be obtained only for {\em some} pairs of species.

In \cite{dre3}, we already
introduced and explored related mathematical questions.
We asked when 
knowing just some of the leaf-to-leaf distances 
is sufficient to uniquely determine -- or, as we say, 
`lasso'  -- the topology of the tree, or its edge weights, or both. In this  paper,
 we turn our attention to 
a fixed (un-weighted) tree $T$ and the set of minimal subsets $\cl$ of $\binom{X}{2}$ for which the leaf-to-leaf distances between all $x,y\in X$ with $\{x,y\}\in \cl$ relative to some edge-weighting $\omega$ of $T$ 
suffice to determine all the other distances relative to $\omega$ and, thus, the edge-weighting $\omega$. Indeed, these  subsets form the bases of the matroid $\rM(T)$ that will be studied here. 

We begin by recalling some basic definitions and some relevant terminology from \cite{dre3} on trees, lassos, and associated concepts (readers unfamiliar with basic matroid theory may wish to consult \cite{wel} -- though even Wikipedia may suffice).  We then define $\rM(T)$ and describe some of its basic properties before presenting our 
main results. Finally, we provide a number of remarks, observations, and questions for possible further study.

\section{ Some terminology and basic facts}
\label{cp-t}

We will assume throughout that $X$ is a finite set of cardinality
$n \geq 3$ and, for any $2$ elements $x,y\in X$, we will usually write just
$xy$ instead of $\{x,y\}$, and we will refer to any such set as a `cord' whenever $x\neq y$ holds. Throughout this paper, we will assume that  $T=(V,E)$ is an $X-$tree, i.e., a
finite tree  with vertex set $V$, leaf set $X\subseteq V$, and edge set $E\subseteq \binom{V}{2}$ that has no vertices of degree $2$. Two $X-$trees $T_1=(V_1,E_1)$ and $T_2=(V_2,E_2)$ are said to be `equivalent' if there exist a bijection $\varphi: V_1 \ra V_2$ with $\varphi(x)=x$ for all $x\in X$ and 
$E_2=\big\{\{\varphi(u),\varphi(v)\}: \{u,v\}\in E_1\}$ in which case we will also write $T_1\simeq T_2.$
In case every interior vertex of an $X-$tree $T$ (that is, every vertex in $V-X$) has degree $3$, $T$ will also be said to be  a `binary' $X-$tree.

Further, given any two vertices $u, v$ of $T$, we denote by $[u,v]$ the set of all vertices on the path $p_T(u,v)$ in $T$ from $u$ to $v$ and by $E(u|v)=E_T(u|v)$ the set of all edges $e$ in $E$ on that path so that $[u,v]=\bigcup_{e\in E(u|v)}e$ always holds.

For each $e\in E$, we denote
by $\omega_e$ the map $E\ra \rr: f\mapsto \delta_{ef}$ (where $\delta$ is, of course, the Kronecker delta function). And  for all $e\in E$ and $xy\in \ch$, we put 
$\delta_{e|xy}:=1$ in case $e\in E(x|y)$ and $\delta_{e|xy}:=0$  otherwise.

Here,  given an $X-$tree $T = (V,E)$,  we will be mainly concerned with the $\rr$-linear map
$$ 
 \rr^E\ra \RR^\ch: \omega\mapsto 
\Big(\omega^T:\ch\ra \rr: xy \mapsto\sum_{e\in E(x|y)}\omega(e) = \sum_{e\in E}\omega(e)\,\delta_{e|xy}\Big) 
$$ 
and the associated $\ch$-labeled family of linear forms
$$ 
\lambda_{xy}=\lambda^T_{xy}: \rr^E\ra \rr:  \omega \mapsto
\omega^T(xy)=
\sum_{e\in E(x|y)}\omega(e) =\sum_{e\in E}\omega(e)\,\delta_{e|xy}
\quad \Big(xy\in\ch\Big).
$$ 
Note that $\lambda_{xy}(\omega_e)=\delta_{e|xy}$ 
holds for all $e\in E$ and all $x,y\in X$, and
$D_\omega(x,y) =\lambda_{xy}(\omega)= \omega^T(xy)$ for all $\omega\in \rr^E$and $x,y\in X$ where $D_\omega=D_{(T,\omega)}$ denotes 
the map
$$
X\times X \rightarrow \RR_{\ge 0}: (x,y)\mapsto D_\omega(x,y) :=\sum_{e\in E(x|y)}\omega(e)
$$
associated to the edge weighting $\omega$ -- a map which in case $\omega$ is a non-negative edge weighting is nothing but the associated (pseudo-)metric on $X$ induced by the edge weighted tree $T=(T,\omega)$ much studied in phylogenetic analysis.

Recall also that, given an arbitrary metric $D:X\times X \ra \RR_{\ge 0}$ defined on $X$, 
\begin{itemize}
\item 
the metric $D$ is dubbed a `tree metric' if it is of the form $D_{(T,\omega)}$ for some $X-$tree $T=(V,E)$ and some non-negative edge weighting $\omega:E \ra \RR_{\ge 0}$ of $T$
\item 
which, in turn, holds if and only if  $D$ satisfies
 the well-known `four-point condition' stating that, for all $a,b,c,d$ in $X$, the larger two of the three distance sums $D(a,b)+ D(c,d), \, D(a,c) + D(b,d), \,D(a,d) + D(b,c)$ coincide,
\item 
that, in this case, one can actually always find an $X-$tree $T$ and an edge weighting 
$\omega$ of $T$ with $D=D_{(T,\omega)}$ such that $\omega$ is strictly positive on all interior edges in which case $\omega$ is called a `proper' edge weighting of $T$, 
\item 
any such pair $(T,\omega)$ is  -- up to canonical isomorphism -- uniquely determined by $D$, 
\item 
and one does not even need  to know the values of $D$ for all cords $xy$ in $\ch$ in order to determine all the other distances and, thus, the edge-weighting $\omega$ in this case.
\end{itemize}

In this note, we continue our investigation of those subsets $\cl$ of $\ch$ for which -- given the $X-$tree 
$T$
-- already the restriction $\omega^T|_\cl$ of the map $\omega^T$ to $\cl$ suffices to determine -- or `lasso' -- the edge weighting $\omega$ of $T$ that we began in \cite{dre3}.
To this end, we denote, for any subset $\cl$ of $\ch$, 

\noindent
\,\, --  by $ \langle \cl \,\rangle= \langle \cl \,\rangle^T$ the $\rr$-linear subspace of the dual vector space 
$\widehat{\RR^E}:=\Hom_\rr(\RR^E,\rr)$ 
of the space $\rr^E$ generated by the maps $\lambda_{xy}$ with $xy \in \cl$,

\noindent
\,\, --  by $\rk(\cl)=\rk^T(\cl):=\dim_\rr \langle \cl \,\rangle$ the dimension of $\langle \cl \,\rangle$, and

\noindent
\,\,-- by $\Gamma(\cl):=(X,\cl)$ the graph with vertex set $X$ and edge set $\cl$.

\noindent
Following the conventions introduced in \cite{dre3}, 

-- we will refer to a subset $\cl$ of $\ch$ as being `connected', `disconnected' or `bipartite' etc. whenever the graph $\Gamma(\cl)$ is connected, disconnected, or bipartite and so on, 

-- a connected component of $\Gamma(\cl)$ will also be called a connected component of $\cl$, 

-- and given any two  subsets $A,B$ of $X$, the subset 
$\{ab: a \in A, b \in B\}$ of $\ch$ will be denoted by $A\vee B$ so that a subset $\cl$ of $\ch$ is bipartite if and only if there exist two disjoint subsets $A,B$ of $X$ with $\cl \subseteq A\vee B$.

Further, a subset $\cl$ of $\ch$ will be called 

-- an `edge-weight lasso' for $T$ if the implication $``\omega_1^T|_\cl=\omega_2^T|_\cl\,\,\,\Ra\,\,\, \omega_1= \omega_2$'' holds for any two proper edge-weightings $\omega_1, \omega_2: E \ra \rr_{>0}$ of $T$,  

-- a `topological lasso' for $T$ if the implication $``\omega_1^{T'}|_\cl=\omega_2^T|_\cl\,\,\,\Ra\,\,\, T'\simeq T \,$'' holds for any $X-$tree $T'$ and  any proper edge-weightings $\omega_1$ of $T'$ and $\omega_2$ of $T$, respectively, and 

-- a `strong lasso' for $T$ if it is simultaneously an edge-weight and a topological lasso for $T$.\\

Next, recall (see\,e.g.\,\cite{ox,wel}) that 
an `abstract' matroid $\rM$ with a ground set, say,  $M$ can be defined in terms of its `rank function' 
$\rk_\rM:\cp(M)\ra \nn_0$ (and with $\cp(M)$ denoting the power set of $M$) as well as by the collection $\rI_\rM=\{I\subseteq M:\rk_\rM(I)=|I|\}$ of its `independent sets', the collection $\rG_\rM=\{L\subseteq M:\rk_\rM(L)=\rk_\rM(M)\}$ of its 'generating sets', the collection $\rB_\rM$ of its `bases', i.e., the maximal sets in 
$\rI_\rM$ or, just as well, the minimal sets in $\rG_\rM$, the collection $\rC_\rM$ of its `circuits', i.e., the minimal sets in $\cp(M)-\rI_\rM$, as well as the `closure operator'  
$[ \dots]_\rM: \cp(M)\ra \cp(M):
L \mapsto [L]_\rM:=\{m \in M: \rk_\rM(L)=\rk_\rM(L\cup\{m\})\}$ 
associated to $\rM$. 

Here, given any $X-$tree $T$, we want to investigate the matroid $\rM(T)$ with ground set $M:=\ch$ associated to $T$ whose rank function $\rk_{\rM(T)}$ is the map 
$\rk^T:\cp(\ch)\ra \nn_0$ defined just above, i.e.,  the matroid that is `represented' (over $\rr$, again see\,e.g.\,\cite{ox,wel}) by the map 
$$
\lambda^T:\ch\ra \widehat{\RR^E}:xy\mapsto \lambda_{xy}.
$$ 
We will denote by $\rI(T):=\rI_{\rM(T)}$ its collection of independent sets, by $\rG(T):=\rG_{\rM(T)}$ its collection of generating sets, 
by $\rB(T):=\rB_{\rM(T)}$ its collection of bases, by $\rC(T):=\rC_{\rM(T)}$ its collection of `circuits' and, given any subset $\cl$ of 
$\ch$, we denote by $[ \cl ]^T :=[\cl]_{\rM(T)}$ the `($T$-)closure' of $\cl$ relative to $\rM(T)$.

\smallskip
\noindent
It was noted already in \cite[Theorem 1]{dre3} that a subset $\cl$ of $\ch$ is an edge-weight lasso for an $X-$tree $T = (V,E)$
if and only if the implication $``\omega_1^T|_\cl=\omega_2^T|_\cl\,\,\,\Ra\,\,\, \omega_1= \omega_2$'' does not only hold for any two proper edge weightings $\omega_1, \omega_2$ of $T$, but for any two maps
$\omega_1, \omega_2\in \RR^E$ and, hence, if and only if 
 $\langle \cl \,\rangle$ coincides with  $\widehat{\RR^E}$ or, using the terminology introduced above, if and only if $\cl \in \rG(T)$  or, just as well, $\rk^T(\cl) =\rk^T(\ch) =|E|$ holds. In particular, an edge-weight lasso $\cl$ for $T$ is a 
`tight' edge-weight lasso for $T$, i.e, a minimal subset of $\ch$ that is an edge-weight lasso for $T$, if and only if its cardinality coincides with $|E|$ if and only if it is a basis of $\rM(T)$, that is, 
$\cl \in \rB(T)$ holds.

\medskip
\noindent

\begin{figure}[h]
 \begin{center}
\resizebox{12cm}{!}{
{
\includegraphics{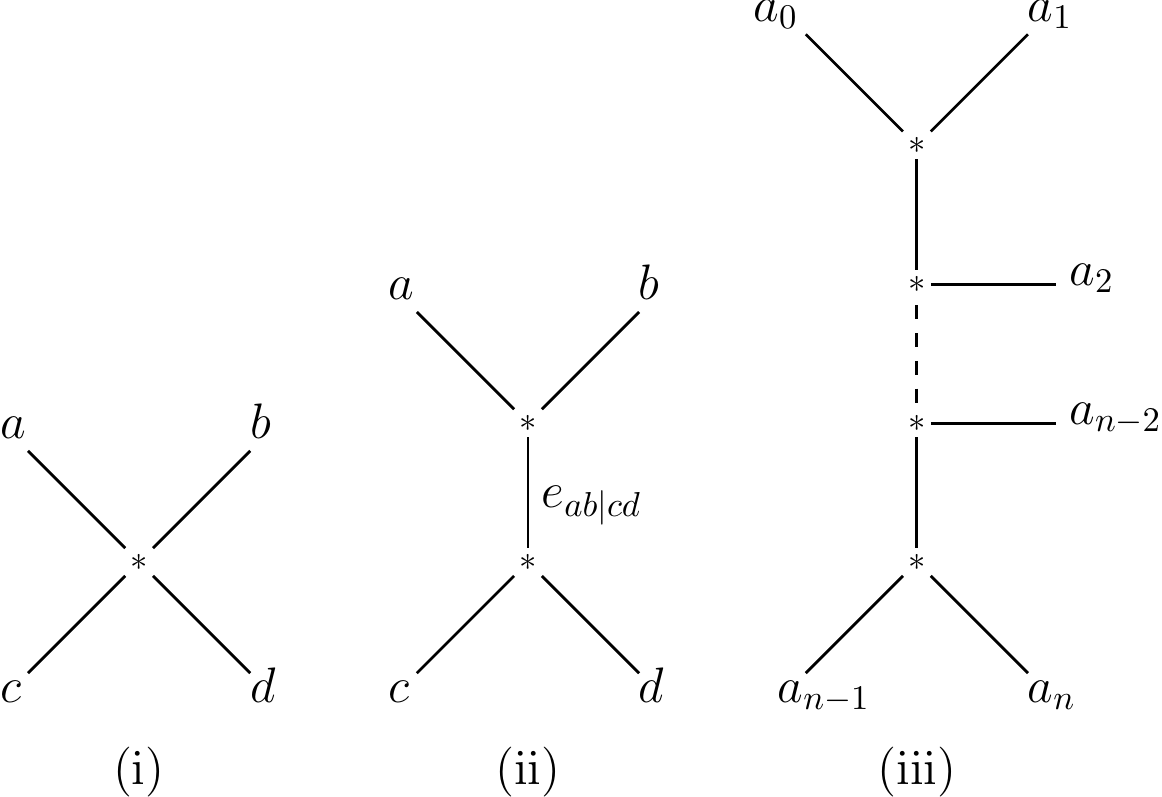}}
}
\caption{\label{figure1}}{(i) A star tree with leaf set $X_4:=\{a,b,c,d\}$; (ii) A binary $X_4-$tree -- up to equivalence, there are two more binary $X_4-$trees ; (iii) a `caterpillar' $X_n-$tree for $X_n:=\{ a_1, a_2, \ldots, a_{n-1}, a_n\}$.}
\end{center}
\end{figure}

Particular types of $X-$trees that will play an important role in this paper are shown in Figure\,\ref{figure1}. They comprise (i) the `star trees', i.e., 
$X-$trees that have
 just one interior vertex and, hence, are equivalent to the tree $T^\star(X):=\big(V^\star(X), E^\star(X)\big)$ with leaf set $X$, vertex set $V^\star(X):=X \dot\cup \{\star\}$ and edge set $E^\star(X):=\big\{\star x:x\in X\big\}$
where `$\star$' denotes just some arbitrary, but fixed element not in $X$;  (ii) `quartet trees', i.e., binary $X-$trees that have four leaves (with $T_{ab|cd}$ denoting the quartet tree with leaf set $\{a,b,c,d\}$ whose central edge that will also be denoted by $e_{ab|cd}$ separates the leaves $a,b$ from $c,d$),
and (iii) `caterpillar trees', i.e. binary $X-$trees $T=(V,E)$ 
containing two interior vertices $u,v\in V-X$ 
with $[u,v]=V-X$.
%



\section{Star trees}
 
For the simplest type of $X-$tree, i.e., the star tree $T^\star:=T^\star(X)$ with leaf set $X$ (cf.\,Figure\,\ref{figure1}\,(i)\big), the associated matroid $\rB(T^\star)$ is well known: It is easily seen to exactly coincide with the `biased matroid' of the complete signed graph $(X,\ch)$ with vertex set $X$ all of whose edges have sign $-1$. 
In consequence (see e.g. \cite[Section 6.10]{ox}
and the references therein to Zaslavsky's papers on signed graphic matroids), the following
results are known to hold:

\begin{proposition}
\label{starprop}
Given a finite set $X$ of cardinality $n\ge 3$, the following holds for the matroid $\rB(T^\star)$ associated to the star tree $T^\star:=T^\star(X)$ with leaf set $X$:
\begin{itemize}
\item[\rm (i)]
The collection $\rG(T^\star)$ of all edge-weight lassos for $T^\star$ coincides with the collection of all `strongly non-bipartite' subsets $\cl$ of $\ch$, i.e, all subsets $\cl$ of $\ch$ for which none of the connected components of $\cl$ is bipartite.
\item[\rm (ii)]
The collection $\rB(T^\star)$ of all tight edge-weight lassos for $T^\star$ coincides with 
the collection of all minimal strongly non-bipartite subsets $\cl$ of $\ch$, i.e, all subsets $\cl$ of $\ch$ for which each connected component of $\cl$ contains exactly one circle\footnote{In our context, we adopt the convention of calling a graph (and, hence, also every subgraph of a graph) a `circle' if it is connected and every vertex in that graph has degree $2$.}
and the length of this circle has odd parity.
\item[\rm (iii)]
The collection $\rI(T^\star)$ of all independent subsets of $\rM(T^\star)$ coincides with 
the collection of all subsets $\cl$ of $\ch$ for which each connected component of $\cl$  is either a tree or contains
exactly one circle and the length of this circle has odd parity. 
\item[\rm (iv)]
The collection $\rC(T^\star)$ of all circuits of $\rM(T^\star)$ coincides with 
the collection of all subsets $\cl$ of $\ch$ that either form a circle of even length or 
a pair of circles of odd length together with a connecting simple path, such that the two circles are either disjoint (then the connecting path has one end in common with each circle and is otherwise disjoint from both) or share just a single common vertex (in this case the connecting path is that single vertex).
\item[\rm (v)]
The co-rank $n-\rk^{T^\star}(\cl)$ of a subset $\cl$ of $\ch$ relative to  $\rM(T^\star)$ coincides with the number of non-bipartite connected component of $\cl$.
\item[\rm (vi)] The closure $[\cl]^{T^\star}$ of a subset $\cl$ of $\ch$ relative to  $\rM(T^\star)$ coincides with the union of 
{\bf (a)} the edge set of the complete graph whose vertex set is the union of the vertex sets of all non-bipartite connected components of $\cl$ and
{\bf (b)} all subsets of the form $A\vee B$ for which some bipartite connected component $\cl'$ of $\cl$ with $\cl'\subseteq A\vee B$, 
$A \cup B= \bigcup_{xy\in \cl'}\{x,y\}$,  and $A \cap B =\0$ exists.
\end{itemize}
\end{proposition}

\section{A recursive approach for computing $\rB(T)$}

Every $X-$tree can be reduced by a sequence of edge contractions to a star tree (one may even insist that at each stage, one of the two subtrees 
incident with the edge being contracted has only one non-leaf vertex, though we do not require this here).  Thus, Proposition~\ref {starprop} can be used as basis for a recursive description of the matroid associated with any $X-$tree,  provided that one can describe, for any $X-$tree $T$, how to obtain
 $\rB(T)$ from $\rB(T/f)$  
 where $f$ is any interior edge of $T$, and $T/f$ is the $X-$tree obtained from $T$ by collapsing edge $f$.
 We provide such a description shortly,  in Proposition~\ref{recursive}, using the following lemma.
 
\begin{lemma}\label{edge-collapsing}
Given any $X-$tree $T=(V,E)$, any subset $F$ of the set of interior edges of $T$, any map $\lambda \in \widehat{\rr^E}$, and any map 
$\rho:\ch \rightarrow \rr$ with 
$\lambda =\sum_{xy\in\ch}\rho(xy)\lambda^T_{xy}$, let 
$T/F$ denote the $X-$tree obtained by collapsing all edges in $F$, and let 
$\lambda|_{E-F}$ denote 
the restriction of $\lambda$ to the space $\rr^{E-F}$ relative to the canonical embedding $\rr^{E-F}\ra\rr^{E}:\omega\mapsto \omega_{(F\ra 0)}$ defined by extending each map 
$\omega\in \rr^{E-F}$ to the map $\omega_{(F \ra 0)}$ by putting $\omega_{(F\ra 0)}(e):=0$ for all 
$e\in F$. Then, one has $\lambda|_{E-F}=\sum_{xy\in\ch}\rho(xy)\lambda^{T/F}_{xy}$, i.e., one has $\lambda(\omega)=\sum_{xy\in\ch}\rho(xy)\lambda^{T/F}_{xy}(\omega|_{E-F})$ for all maps $\omega\in \rr^E$ with
$\omega(f)=0$ for all $f\in F$. 

In particular, given any edge-weight lasso $\cl$ for $T$,  $\cl$ is also an edge-weight lasso 
for the $X-$tree $T/F$.
More generally, 
 \begin{equation}\label{rank/F}
\rk^T(\cl)=\rk^{T/F}(\cl) + \dim\{\lambda\in \langle \cl \,\rangle^T:\lambda(e)=0 \text { for all } e \in E-F\}
\le \rk^{T/F}(\cl) + |F|
\end{equation}
holds for every subset $\cl$ of $\ch$ and any subset $F$ of the set of interior edges of $T$.
\end{lemma}

\pf
The first part follows directly from the definitions and implies that $\lambda_{xy}^{T/F}=\lambda_{xy}^T|_{E-F}$ holds 
for all $xy\in \ch$. In particular, as the map $ \widehat{\rr^{E}}\ra \widehat{\rr^{E-F}}: \lambda \mapsto \lambda|_{E-F}$ is surjective,
the maps 
$\lambda_{xy}^{T/F} \,\,(xy \in \ch)$ must generate $\widehat{\rr^{E-F}}$ whenever the maps 
$\lambda_{xy}^T \,\,(xy \in \ch)$ generate $\widehat{\rr^E}$ while, more generally, they generate a space whose dimension coincides with the difference of $\rk^T(\cl)$ and the dimension of the kernel of the map $\langle \cl \,\rangle^T\ra \widehat{\rr^{E-F}}: \lambda \mapsto  \lambda|_{E-F}$.
\epf

\begin{proposition}
\label{recursive}
Given an $X-$tree 
$T$,
an interior edge $f$ of $T$, a pair $xy \in \ch$, and a basis $B$ of $\rM(T/f)$, let 
$\rho_{xy}\in \rr^{B}$
  denote the unique map in  
 $\rr^{ B}$  
with 
 $\lambda^{T/f}_{xy} = \sum_{b \in B} \rho_{xy}(b) \lambda^{T/f}_b$. 
Then, $\rB(T)$ coincides with the set 
$$
\rB_{f}:= \Big\{ \{xy\} \cup B:  xy\in \ch, B \in \rB(T/f),\mbox{ and }
 \sum_{b \in B} \rho_{xy}(b) \,\delta_{fb}\neq \delta_{f|xy}\Big\}.
 $$
\end{proposition}
\pf
By Lemma~\ref{edge-collapsing}, there exists, for each $B' \in \rB(T)$, some $b' \in B'$ with $B'-b' \in \rB(T/f)$.  Thus, each element of $\rB(T)$ is of the form $B \cup \{xy\}$ for some
$xy \in \ch$ and some $B \in \rB(T/f)$. So, 
$\rB(T)\subseteq \big\{ \{xy\} \cup B:  xy\in \ch, B \in \rB(T/f)\big\}$ must clearly hold. Now suppose that $xy\in \ch$ and $B \in \rB(T/f)$ holds.
Then, denoting by $\Lambda_{B,xy}$ the space of all maps $\lambda \in   \langle  \{xy\} \cup B \rangle^T$ with $\lambda|_{E-\{f\}}=0$, it follows from (\ref{rank/F}) that 
$
\{xy\} \cup B\in \rB(T) \iff \dim \Lambda_{B,xy}=1
$
holds while, by construction, we have 
$ \dim \Lambda_{B,xy}=1$ if and only if there exists some non-zero map in $\Lambda_{B,xy}$.

However, given any map $\rho\in \rr^B$ and any real number $c$, it follows from the fact that, by definition,
$\lambda^{T}_{zz'}|_{E-\{f\}}$ coincides with $\lambda^{T/f}_{zz'}$ for all $zz'\in \ch$, one has 
 $\lambda_{c,\rho}:=-c \lambda^{T}_{xy}+\sum_{b \in B} \rho(b) \lambda^{T}_b\in\Lambda_{B,xy}$ if and only if 
$-c \lambda^{T/f}_{xy} +\sum_{b \in B} \rho(b) \lambda^{T/f}_b$ vanishes and, hence, if and only if $\rho(b) =c \rho_{xy}(b) $ holds for all $b\in B$. Thus, one has $\{xy\} \cup B\in \rB(T)$ if and only if one has  $\lambda_{c,\rho}(\omega_f)\neq 0$ and, hence, if and only if $\sum_{b \in B} \rho_{xy}(b) \,\delta_{fb}\neq \delta_{f|xy}$ holds as claimed. \epf

\noindent
{\bf Remark:}

Similarly, suppose that $T=(V,E)$ is an $X-$tree and that $U\subseteq V$ is a $T-$core as defined in \cite[Section 5]{dre3}, i.e., a non-empty subset of $V$ for which  the induced subgraph $T_U:=(U,E_U:=\{e \in E: e \subseteq U\})$ of $T$ with vertex set $U$ is connected (and, hence, a tree) and the degree $\deg_{T_U}(v)$ of any vertex $v$ in $T_U$ is either $1$ or coincides with the degree $\deg_T(v)$ of $v$ in $T$. Then,  the rank $\rk^T(\cl)$ of a subset $\cl$ of $\ch$ relative to $T$ and the  rank 
$\rk^{T_U}(\cl_U)$  of the corresponding subset $\cl_U$ of $\binom{X_U}{2}$ relative to the $X_U-$tree $T_U$ are easily seen to be related by the inequality
$$
\rk^T(\cl) \le \rk^{T_U}(\cl_U) + |E-E_U|.
$$
This fact can be used to prove  \cite[Theorem 5]{dre3} in the same way Lemma \ref{edge-collapsing} has been used above to establish Proposition \ref{recursive}.

\subsection{An example}

\begin{figure}[h] \begin{center}

\resizebox{12cm}{!}{
{
\includegraphics{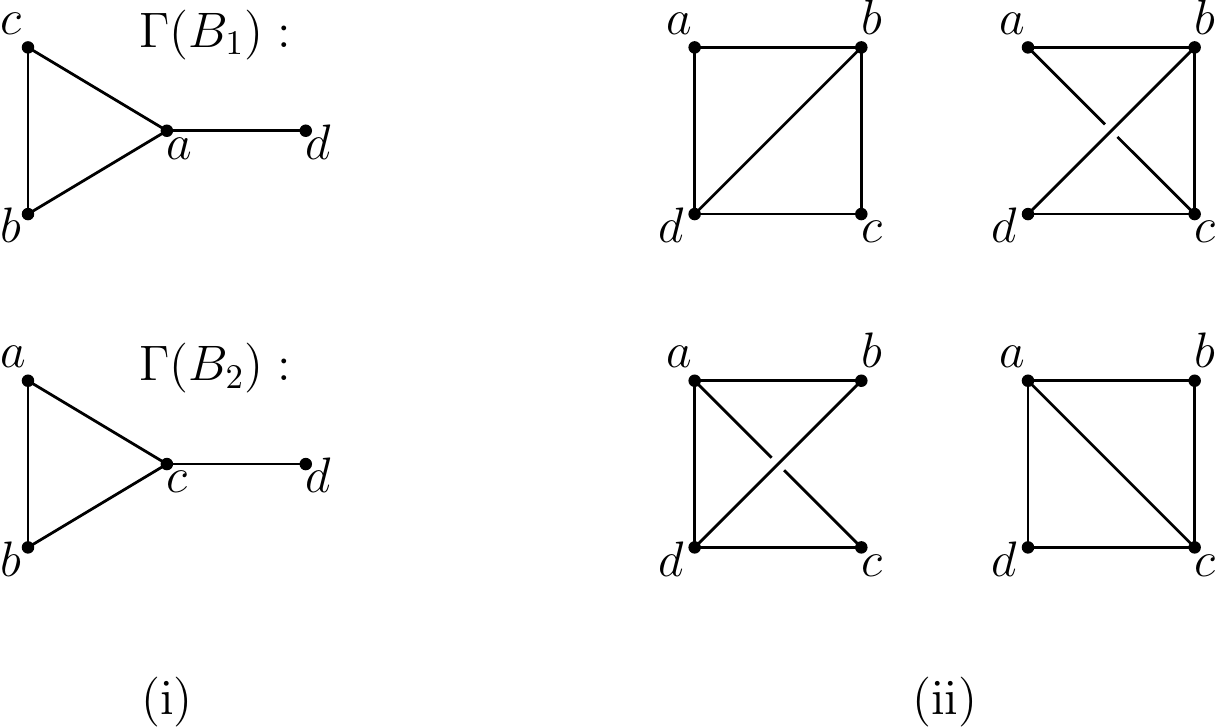}}
}
\caption{\label{figure2}}{(i): Two graphs representing two of the twelve tight edge-weight lassos of $T^\star(\{a,b,c,d\})$, one from each of 
the two orbits of such lassos relative to the $8$-element symmetry group of $T_{ab|cd}$. 

(ii): The four graphs associated to the four bases in 
 $\rB(T_{ab|cd})$.}

\end{center}
\end{figure}

To illustrate Proposition~\ref{recursive}, consider -- for $X:=\{a,b,c,d\}$ -- the quartet $X-$tree $T:=T_{ab|cd}$ shown in 
Figure\,\ref{figure1}\,(ii).  In this case, there is -- up to scaling -- only one linear relation between the six maps 
$\lambda_{xy}^T \,\,(xy \in \ch)$, viz. the relation 
$$
\lambda_{ac}^T +\lambda_{bd}^T=\lambda_{ad}^T +\lambda_{bc}^T: E^\rr\ra \rr: \omega \mapsto 2 \omega(e_{ab|cd}) + \sum_{e\in E, e \neq e_{ab|cd}}  \omega(e).
$$ 
Thus, $\rB(T)$ consists of the four $5$-subsets $\cl$ of $\ch$ 
that do not contain exactly one of the four cords $ac,ad,bc,bd$ -- or, equivalently, with $|\cl\cap \{ac,ad,bc,bd\}|=3$ -- and, hence, the four subsets $\cl$ of $\ch$ 
whose graphs $\Gamma(\cl)$ are shown in 
Figure\,\ref{figure2}(ii).  Clearly, if $f$ coincides with the unique interior edge 
 of $T_{ab|cd}$ \big(i.e., the edge denoted by
$e_{ab|cd}$ in 
Figure\,\ref{figure1}\,(ii)\big), $T/f$ is equivalent to the star tree 
$T^\star:=T^\star(X)$, also shown in Figure\,\ref{figure1}\,(i), and the graphs $\Gamma(\cl)$ corresponding to the bases $\cl$ in $\rB(T^\star)$, being minimal strongly non-bipartite graphs with vertex set $X$, must consist of one triangle (for which there are $4$ possibilities) to which the remaining element in $X$ is appended by a single edge  (for which there are $3$ possibilities). So,  
$\rB(T^\star)$ consists of $12$ bases that form two orbits relative to the symmetry group of $T$ representatives of which are the bases $B_1:=\{ab,bc,ca,da\}$ and $B_2:=\{ab,bc,ca,dc\}$ shown in Figure\,\ref{figure2}\,(i). For the two cords $db,dc\in \ch - B_1$, 
we have -- putting $f:=e_{ab|cd}$ --
$$
\lambda^{T/f}_{db}=\lambda^{T/f}_{da}-\lambda^{T/f}_{ac}+\lambda^{T/f}_{cb}
\text{ \quad and \quad} 
\lambda^{T/f}_{dc}=\lambda^{T/f}_{da}-\lambda^{T/f}_{ab}+\lambda^{T/f}_{bc}
$$ 
while
$$
\lambda^{T}_{db}(\omega_f)= \lambda^{T}_{da}(\omega_f)-\lambda^{T}_{ac}(\omega_f)+\lambda^{T}_{cb}(\omega_f) = 1
$$
and
$$
\lambda^{T}_{dc}(\omega_f)=0\neq \lambda^{T}_{da}(\omega_f)-
  \lambda^{T}_{ab}(\omega_f)+
 \lambda^{T}_{bc}(\omega_f) =
 2
$$
holds implying that, to bases of type $B_1$, we can add cords of type $dc$, but not cords of type $db$.\\

And for the two cords $da,db\in \ch - B_2$, we have 
$$
\lambda^{T/f}_{da}=\lambda^{T/f}_{dc}-\lambda^{T/f}_{cb}+\lambda^{T/f}_{ba}
\text{ \quad and \quad} 
\lambda^{T/f}_{db}=\lambda^{T/f}_{dc}-\lambda^{T/f}_{ca}+\lambda^{T/f}_{ab}
$$ 
as well as
$$
1=\lambda^{T}_{da}(\omega_f)\neq \lambda^{T}_{dc}(\omega_f)-\lambda^{T}_{cb}(\omega_f)+\lambda^{T}_{ba}(\omega_f) = -1
$$
and
$$
1=\lambda^{T}_{db}(\omega_f)
\neq \lambda^{T}_{dc}(\omega_f)-
  \lambda^{T}_{ca}(\omega_f)+
 \lambda^{T}_{ab}(\omega_f) =-1
$$
holds implying that, to bases of type $B_2$, we can add either one of the two missing cords.
Obviously, this fully
corroborates our previous assertion about $\rB(T_{ab|cd})$.\\

{

 \section{Some particular cases}
 
  \subsection{Pointed $x$-covers of binary $X-$trees $T$ that are bases of $\rM(T)$}
 
 When 
$T$
is a binary $X-$tree, some particular bases in $\rB(T)$ are easily described: Select any element $x \in X$ and, for each one of the 
 $n-2$ interior vertices $v$ of $T$, consider
 the three components of the graph obtained from $T$ by deleting $v$.  Select an element of $X$ from each of the two components that do not contain $x$, and denote this pair by $y_v=y_v(x), z_v=z_v(x)$.  Put
 $$
 \cl= \{ax: a \in X-\{x\}\} \cup \{y_vz_v: v \in V-X\},
 $$ 
 and let $P_x(T)$ denote 
 the collection of 
 subsets of $\ch$ that can be generated in this way (by the various choices of $y_v$ and $z_v$ as $v$ varies). 

 For example, considering again the quartet $X-$tree $T:=T_{ab|cd}$ with its two interior vertices $u$ and $v$ as shown in Figure\,\ref{figure1}\,(ii), we may choose $x:=d$,  $y_u=y_v:=a, z_u:=b, z_v=c$ and obtain the lasso 
$$
\cl= \{ad,bd,cd\} \cup \{ab,ac\}
$$
 as an element of  $P_d(T_{ab|cd})$.
 
Clearly, $P_x(T)$ is a subset of $\rB(T)$ for each $x \in X$, since the elements of $P_x(T)$ correspond precisely to the so-called  `pointed $x-$covers' of $T$ of cardinality $2n-3$ and, by
 Theorem 7 of \cite{dre3}, any pointed $x-$cover $\cl$ of a 
 binary $X-$tree is not only an edge-weight, but a strong lasso for that tree.


 We note also that, given two distinct elements  $x_1,x_2$ in 
$X$, a subset $\cl$ of $\ch$ cannot simultaneously be a pointed $x_1$-cover in $P_{x_1}(T)$ and a pointed $x_2$-cover in $P_{x_2}(T)$ unless $T$ is a caterpillar tree with $x_1$ and $x_2$ at opposite `ends' of the tree: Indeed, if there  exists some
$\cl\in P_{x_1}(T)\cap P_{x_2}(T)$, 
we must have
 $\{y_v(x_1)z_v(x_1): v \in V-X\} =\big\{x_2a:a\in X-\{x_1,x_2\}\big\}$ implying that the path from $x_1$ to $x_2$ in $T$ must pass through every interior vertex of $T$.\medskip 
 

Our next results require two definitions that will also be important later in this paper: 
Recall first that, given an $X-$tree 
$T$ 
and a subset $Y$ of $X$ of cardinality at least $3$, one denotes 
\begin{itemize}
\item 
by $T|_Y$  
the $Y-$tree obtained from the minimal subtree of $T$ that connects the leaves in $Y$ by suppressing any resulting vertices of degree 2 (see e.g.\,\cite[Section 2.3]{dre3}), 
\item 
by $V|_Y$ and $E|_Y$ its vertex and 
edge set, respectively,
\item 
and, given in addition any edge weighting $\omega$ of $T$, one denotes by $\omega|_Y$ the `induced'  edge weighting of $T|_Y$, i.e., the edge weighting that maps any edge 
$\{u,v\}\in E|_Y$ onto the sum $\sum_{e\in E(u|v)} \omega(e)$, yielding a surjective $\rr$-linear map 
$\res_Y:\rr^E\ra \rr^{E|_Y}:\omega\mapsto \omega|_Y$ such that 
$\lambda^T_{yy'}(\omega)=\omega^T(yy')= (\omega|_Y)^{T|_Y}(yy')=\lambda^{T|_Y}_{yy'}\big(\res_Y(\omega)\big)$ holds for all $\omega \in\rr^E$ and all $yy'\in \binom{Y}{2}$.
\end{itemize}
It follows that the map 
$\lambda^T_{yy'}:\rr^E\ra \rr$ coincides, for all $yy'\in \binom{Y}{2}$, 
with the map $\lambda^{T|_Y}_{yy'} \!\circ \res_Y$, the composition of the maps $\res_Y$ and $\lambda^{T|_Y}_{yy'}$.

So, denoting by 
 $\widehat{\res_Y}$ the dual -- and necessarily injective -- map $\widehat{\rr^{E|_Y}}\ra \widehat{\rr^E}: \lambda \mapsto \lambda \circ \res_Y$ of the map $\res_Y$, we have 
also 
$\langle \cl \,\rangle^T= \widehat{\res_Y}(\langle \cl \,\rangle^{T|_Y})$ and $\rk^T(\cl)=\rk^{T|_Y}(\cl)$ for every subset $\cl$ of $\binom{Y}{2}$.
In consequence,
we must also have 
\begin{equation}\label{res_Y}
\rk^T(\cl)=\rk^{T|_Y}(\cl)
\end{equation}
for every subset $Y$ of $X$ of cardinality at least $3$ and every subset $\cl$ of $\binom{Y}{2}$, implying also that 
 every circuit $\cl\subseteq \binom{Y}{2}$ of $\rM(T|_Y)$ must also be a
 circuit of $\rM(T)$, i.e., we have 
$\rC(T|_Y)\subseteq \rC(T)$ for every such subset $Y$ of $X$ of cardinality at least $3$.\\

Further,  denoting -- for every $x\in X$ -- by 
 $e_x\in E$ the unique  pendant edge of $T$ containing $x$,
we say that a $2$-subset $ab$ of $X$ forms (or `is') a `$T-$cherry' 
if the two edges $e_a,e_b\in E$ share a vertex, and  $ab$ is said to form a `proper $T-$cherry' if this vertex has degree 3.
Note that a 
$2$-subset $ab$ of $X$ forms a proper $T-$cherry if and only if $T|_{\{a,b,c,d\}}\simeq T_{ab|cd}$ holds for any two distinct elements $c,d$ in $X-\{a,b\}$ (if any). Note also that, in a binary $X-$tree $T$, every $T-$cherry is a proper $T-$cherry. In addition, such a tree is a caterpillar tree if and only if $n=3$ holds or its leaf set contains exactly two distinct $T-$cherries.  We claim:
\begin{proposition}\label{co-loops}
For an $X-$tree $T$,  a cord $ab \in \ch$ is a `co-loop' of $\rM(T)$, i.e.  it is contained in every edge-weight lasso for $T$, if and only if $ab$ is a proper $T-$cherry.
\end{proposition}
\pf
 If $n= 3$ holds, the set $\ch$ is the only basis of $\rM(T)$ while,  if  $n\ge 4$ holds and
$ab$ is a proper $T-$cherry, the cord $ab$ must be contained in every edge-weight lasso for $T$ in view of \cite[Corollary 1]{dre3}.
Conversely, if $ab$ does 
not form a proper $T-$cherry, there must exist two distinct elements $c,d$ in $X-\{a,b\}$ such that 
$\lambda_{ab}^T+\lambda_{cd}^T=\lambda_{ac}^T+\lambda_{bd}^T$ holds, implying that
$\lambda_{ab}^T$ cannot be a co-loop of $\rM(T)$. \epf

\section{Main results}

\subsection{$\rM(T)$ determines $T$ up to equivalence}
We begin this section by showing that the matroid associated with an $X-$tree determines that $X-$tree up to equivalence:
\begin{theorem}
\label{mthm}
One has ``$\rM(T_1)=\rM(T_2) \Longleftrightarrow T_1\simeq T_2$'' for 
 any two $X-$trees $T_1$ and $T_2$.
\end{theorem}
\pf
We first note that, if $T$ is any $X-$tree and $Y=\{a,b,c,d\}$ is a $4$-subset of $X$, then we have 
$T|_Y \simeq T_{ab|cd}$ if and only if there exists at  least one basis $B$ of $\rM(T)$ containing the set $\cl_{abcd}:=\{ab,bc, cd, da\}$: 
Indeed, if $T|_Y \simeq T_{ab|cd}$ holds, the four maps 
$\lambda^{T|_Y}_{xy} \,\,(xy\in \cl_{abcd})$
and, hence, also the corresponding four maps 
$\lambda^T_{xy}\,\,(xy\in \cl_{abcd})$ 
are linearly independent. So, by the matroid augmentation property of independent sets, there exists some  $B \in \rM(T)$ containing these four cords.
 Conversely, if $T|_Y\not\simeq T_{ab|cd}$ and, therefore, also 
$ \lambda^T_{ab}+\lambda^T_{cd} =\lambda^T_{ad}+ \lambda^T_{bc}$ holds, $\cl_{abcd}$ cannot be part of a basis $B \in \rM(T)$.
It follows that $\rM(T_1)\simeq \rM(T_2)$ implies $\mathbb  Q(T_1) = \mathbb Q(T_2)$ where, for any $X-$tree $T$,  
$\mathbb Q(T)$ is defined by  $\mathbb Q(T):= \{ab|cd:  \{a,b,c,d\} \in \binom{X}{4}, T|_{\{a,b,c,d\}} \simeq T_{ab|cd}\}$. 
However, it has been
observed already by H. Colonius and H. Schultze in \cite{col77,col81}
that $\mathbb Q(T_1) = \mathbb Q(T_2)$ holds for any  two $X-$trees $T_1, T_2$ if and only if one has $T_1 \simeq T_2$ (for a more recent account, see \cite[Theorem 2.7]{dre11}) or \cite[Corollary 6.3.8]{sem}. 
\epf

\subsection{The rank of topological lassos}
Now assume that $n\ge 4$ holds and recall that the following three assertions are -- according to \cite[Theorem 8]{dre3} -- equivalent in this case for any $X-$tree 
$T=(V,E)$
 and any bipartition of $X$ into two disjoint non-empty subsets $A,B$:
\begin{itemize}
\item[{\rm (split-i)}] The subset $A\vee B$ of $\ch$ is
a topological lasso for $T$,

\item[{\rm (split-ii)}] $A\vee B$ is  a `$t-$cover' of $T$ (i.e., given any interior vertex $v$ of $T$ and any two edges $e,e'\in E$ with $v\in e,e'$, there exists some  cord $xy$ in $\cl$ with $e,e'\in E(x|y)$, see \cite[Section 7]{dre11}).
\item[{\rm (split-iii)}] $A \cap \{a,b\}\neq \emptyset \neq B \cap  \{a,b\}$ holds for every 
$T-$cherry $ab$.\footnote{When stating this theorem in \cite{dre3}, we forgot to mention that one needs to assume that $n\ge 4$ holds. Indeed, it is simply wrong  for $n=3$ for obvious trivial reasons as (split-i) holds for all bipartitions $A,B$ of the leaf set $X$ of an $X-$tree with $3$ leaves, but (split-ii) and (split-iii) never holds in this case. Yet, the assumption $n\ge 4$ will always be made here when applying this theorem.}
\end{itemize}
And it was also noted in this context that such bipartitions exist if and only if every $T-$cherry is a proper $T-$cherry and $n\ge 4$ holds.

Here, we want to complement this result as follows:

\begin{theorem}\label{bipartite}
Given any $X-$tree $T=(V,E)$, one has $\rk^T(\cl)\le |E|-1$ for every 
 bipartite subset $\cl$ of $\ch$. Furthermore, the following assertions are equivalent for every such subset $\cl$ of $\ch$:
\begin{itemize}
\item[\rm(i)]
The rank $\rk^T(\cl)$ of $\cl$ coincides with $|E|-1$.

\item[\rm(ii)]  There exists some cord $xy\in \ch$ such that $\cl\cup\{xy\}$ is an edge-weight lasso for $T$.

\item[\rm(iii)]  
$\cl$ is connected and $\cl\cup\{xy\}$ is an edge-weight lasso for $T$ for some cord $xy\in \ch$ if and only if $\cl\cup\{xy\}$ is not bipartite.
    
\item[\rm(iv)] $\cl$ is connected, the closure $[\cl]^T$ of $\cl$ relative to $\rM(T)$  coincides with the edge set of the $($necessarily unique$)$ complete bipartite graph with vertex set $X$ whose edge set contains $\cl$, i.e., $[\cl]^T$ coincides with the set 
$A \vee B$ in case the two subsets  $A,B$ of $X$ form the 
$($necessarily unique$)$ bipartition of $X$ with $\cl\subseteq A \vee B$,
and this set forms a `hyperplane' in $\rM(T)$, i.e., a maximal subset of $\ch$ of rank smaller than $|E|$.
\end{itemize}
\end{theorem}

\pf Assume that $A$ and $B$ are two subsets of  $X$ that form a bipartition of $X$ with $\cl\subseteq A \vee B$, and let $\omega_{A|B}\in \rr^E$ denote the map in $\rr^E$ that maps every interior edge of $T$ onto $0$, every pendant edge $e$ that is incident with some leaf in $A$ onto $1$, and every pendant edge $e$ that is incident with some leaf in $B$ onto $-1$. Clearly,  
\begin{equation}\label{AB}
\lambda_{xy}^T(\omega_{A|B})=\begin{cases}
      \,\,+ 2 & \text{ if } x,y\in A, \\
      \,\,-2& \text{ if } x,y\in B, \\
      \quad 0 &  \text{ otherwise},
\end{cases}
\end{equation}
 holds for every cord $xy\in \ch$.
In particular, one has 
$\lambda_{xy}^T(\omega_{A|B})=0$ for some cord $xy\in \ch$ if and only if
$xy\in A \vee B$ holds. So, standard matroid theory implies that $\rk^T(\cl)\le \rk^T(A \vee B)\le |E|-1$ must hold for every 
subset $\cl$ of a set of the form $A\vee B$ for some bipartition $A,B$ of $X$,
that is,
for every bipartite  subset $\cl$ of $\ch$ -- which is just our first assertion.

(i)\,$\LRa$\,(ii): It follows also from standard matroid theory that the rank $\rk^T(\cl)$ of an arbitrary nonempty subset $\cl$ of $\ch$ with  $\rk^T(\cl)\le |E|-1$ coincides with $|E|-1$ if and only if there exists some cord $xy\in \ch$ such that $\cl\cup\{xy\}$ is an edge-weight lasso for $T$.

(i)\,$\Ra$\,(iii): And standard matroid theory implies also that, if $\cl$ is any subset of $\ch$, one has $\rk^T(\cl)= |E|-1$ if and only if there exists -- up to scaling -- exactly one non-zero map $\omega\in \rr^E$ with $\lambda_{ab}^T(\omega)=0$ for all $ab\in \cl$ and that, in this case, $\cl\cup\{xy\}$ is an edge-weight lasso for $T$ for some cord $xy\in \ch$ if and only if $\lambda_{xy}^T(\omega)\neq 0$ holds for this map $\omega$.

In consequence,  if $\cl\subseteq \ch$ is bipartite and $\rk^T(\cl)= |E|-1$ holds, our observations above imply that there must be a unique bipartition $A,B$ of $X$ with $\cl\subseteq A \vee B$ and that, given any cord $xy\in \ch$,
the union $\cl\cup\{xy\}$ is an edge-weight lasso for $T$ if and only if 
$\lambda_{xy}^T(\omega_{A|B})\neq 0$ and, hence, if and only if $xy\not\in A \vee B$ holds. Furthermore, the fact that there is only one bipartition $A,B$ of $X$ with $\cl\subseteq A \vee B$ implies that $\cl$ must be connected and that, in consequence, $xy\not\in A \vee B$ holds if and only if $\cl\cup\{xy\}$ is not bipartite.

So, we see that  $\cl\cup\{xy\}$ is indeed an edge-weight lasso for $T$ if and only if  $\cl\cup\{xy\}$ is not bipartite, as claimed.

(iii)\,$\Ra$\,(ii): This is trivial as any connected graph with at least three edges can be extended by a single edge to become a non-bipartite graph.

(i,ii,iii)\,$\Ra$\,(iv): If $\cl$ is connected and bipartite, there exists exactly one bipartition $A,B$ of $X$ with $\cl\subseteq A \vee B$ and $\cl\cup\{xy\}$ bipartite for some cord $xy\in \ch$ if and only if $xy\in  A\vee B$ holds. So, if in addition also $|E|-1 =\rk^T(\cl)$ holds, $ [\cl]^T$ must be a hyperplane and we must have
$|E|-1 =\rk^T(\cl)\le  \rk^T(A\vee B)\le  |E|-1$ and, therefore, $\rk^T(\cl)= \rk^T(A\vee B)= |E|-1$ as well as 
\begin{eqnarray*}
 [\cl]^T&=&\{xy\in \ch:  \rk^T(\cl\cup\{xy\})\le  |E|-1\}\\
&=& \{xy\in \ch: \cl\cup\{xy\} \text{ is bipartite}\}
= A\vee B, 
\end{eqnarray*}
as claimed.

(iv)\,$\Ra$\,(i): This is trivial in view of the fact that $ \rk^T(\cl)= \rk^T([\cl]^T)$ holds for every nonempty subset $\cl$ of $\ch$ and the fact that any maximal subset of $\ch$ of rank smaller than $|E|$ must have rank  $|E|-1$. \epf

The above theorem has an interesting application regarding  topological lassos:
\begin{theorem}\label{topological}
{\rm (i)} Given any $X-$tree $T=(V,E)$ and any 
bipartite subset $\cl$ of $\ch$ with $\rk^T(\cl)=|E|-1$,
the hyperplane $[\cl]^T$ is a topological lasso for $T$ if and only if every $T-$cherry is a proper $T-$cherry 
$($i.e., if and only
if there exists at least one bipartition $A',B'$of $X$ such that 
 $A' \vee B'$ is a topological lasso for $T)$ in which case  $[\cl]^T\cup\{xy \}$ must be a 
 strong lasso for $T$ for every
cord $xy \in \ch$ for which  $(X,\cl\cup\{xy \})$ is not bipartite, that is, with $xy\not\in [\cl]^T$.

{\rm (ii)} Conversely, if $\cl$ is any topological lasso for $T$ with $\rk^T(\cl)<|E|$, then 
$\cl$ must be  bipartite, $\rk^T(\cl) = |E|-1$ must hold, every $T-$cherry must be a proper $T-$cherry, and $\cl\cup\{xy \}$ -- and not only $[\cl]^T\cup\{xy \}$ -- must be a 
 strong lasso for $T$ for every
cord $xy \in \ch$ for which  $\cl\cup\{xy \}$ is not bipartite.

\end{theorem}
\pf
(i): In view of the equivalence (i) $\iff $ (iv) of  Theorem \ref{bipartite}, there must exist a (necessarily unique) bipartition of $X$ into two disjoint subsets $A$ and $B$ such that the hyperplane $[\cl]^T$ coincides with the set $A \vee B$.
Furthermore, this set must also be a topological lasso for $T$ if every $T-$cherry is a proper $T-$cherry: Indeed, 
in view of the results from \cite{dre3} quoted above, it suffices to show that, if $ab$ is a proper $T-$cherry and $a\in A$ holds, we must have $b\in B$. Yet, otherwise, we would have $b\in A$ and, therefore, $ab\not\in\cl$ which would, in case $n\ge 4$, allow us to construct yet another non-zero map $\omega_{ab} \in \rr^E$ with $\lambda_{xy}^T(\omega_{ab})=0$ for all cords $xy\in \cl$ that is not a scalar multiple of $\omega_{A|B}$: Indeed, if the two edges $e_a,e_b\in E$ containing $a$ and $b$, respectively, share the vertex $v$, there would exist exactly one further edge $e_{ab}\in E$ with $v\in e_{ab}$, and putting $\omega_{ab}(e_a)=\omega_{ab}(e_b):=1,\, \omega_{ab}(e_{ab}):=-1$, and 
$\omega_{ab}(e):=0$ for all other edges $e\in E$ would yield such a map $\omega_{ab}$, as required. So, 
$[\cl]^T=A \vee B$ must indeed be a topological lasso for $T$, as claimed.
 
(ii) Conversely, assume that $\cl$ is a topological lasso for $T$ of rank less than $|E|$. Then, there must exist a non-zero map $\omega_0\in \rr^E$ with $\lambda_{xy}(\omega_0)=0$ for all $xy\in \cl$. If 
$\omega_0(e_0)\neq 0$ held for some interior edge $e_0\in E$, we could find some proper edge-weighting 
$\omega\in \rr^E_{\ge 0}$ of $T$ with
$\omega(e_0) =  |\omega_0(e_0)|$, 
while
$\omega(e) >  \omega_0(e)$ holds for all edges $e\in E-\{e_0\}$ which -- in turn -- would imply that the map $\omega':=\omega-\sgn\big(\omega_0(e_0)\big)\omega_0$ would be a map in $\rr_{\ge 0}^E$ that would be a non-proper 
 edge-weighting of $T$ for which 
 $D_{(T,\omega)}|_\cl =D_{(T,\omega')}|_\cl$
holds. In view of the last remark in \cite[Subsection 2.2]{dre3}, this would contradict our assumption that $\cl$ is a topological lasso for $T$. \\
  

 Consequently, 
 the support $\supp(\omega_0):=\{e\in E: \omega_0(e)\neq 0\}$ of $\omega_0$ must be contained in the set $\{e_x:x\in X\}$ of pendant edges of $T$, implying that 
 $\lambda_{xy}(\omega_0)=\omega_0(e_x)+ \omega_0(e_y)$ must hold for every cord $xy\in \ch$. Thus, putting 
$A:=\{x\in X: \omega_0(e_x) >0\}$ and 
 $B:=\{y\in X: \omega_0(e_y) <0\}$, we see that ``$x\in A\iff y\in B$'' must hold for every cord $xy\in \ch$ with $\lambda_{xy}(\omega_0)=0$ and, hence, for all $xy\in \cl$.

Thus, as $\cl$ must be connected for every topological lasso $\cl$ for $T$ in view of \cite[Theorem 4]{dre3}, it follows that the pair $A,B$ of subsets of $X$ forms a bipartition of $X$, that $\cl$ must be bipartite relative to this partition, that $A\vee B$ must also be a topological lasso for $T$, and that -- in consequence -- every $T-$cherry must be a proper $T-$cherry and $\omega_0$ must be a (positive) scalar multiple of the map $\omega_{A|B}\in\rr^E$ defined above.
In particular, there can be -- up to a scaling --  only one non-zero map $\omega\in \rr^E$ 
with $\lambda_{xy}(\omega)=0$ for all $xy\in \cl$ implying that the rank of $\cl$ must indeed coincide with $|E|-1$ and that
$\cl\cup \{xy\}$ must, therefore, be a strong lasso for $T$ for every
cord $xy\in \ch$ for which  $\cl\cup \{xy\}$ is not bipartite, i.e. for every
cord $xy\in \ch-A\vee B$.
\epf

\begin{corollary}\label{proper}
Given any $X-$tree 
$T=(V,E)$, the following assertions are equivalent:
\begin{itemize}
\item[$(i)$] There exists a bipartite subset $\cl$ of $\ch$ that is a topological lasso for $T$.
\item[$(ii)$] There exists a topological lasso $\cl$ for $T$ with $\rk^T(\cl)<|E|$.
\item[$(ii')$] There exists a topological lasso  $\cl$ for $T$ with $\rk^T(\cl)=|E|-1$.
\item[$(iii)$] Every $T-$cherry is a proper $T-$cherry.
\end{itemize}
\end{corollary}


A simple example to illustrate Corollary \ref{proper} is presented in Figure \ref{figure3}:

\begin{figure}[h]
 \begin{center}
\resizebox{12cm}{!}{
{
\includegraphics{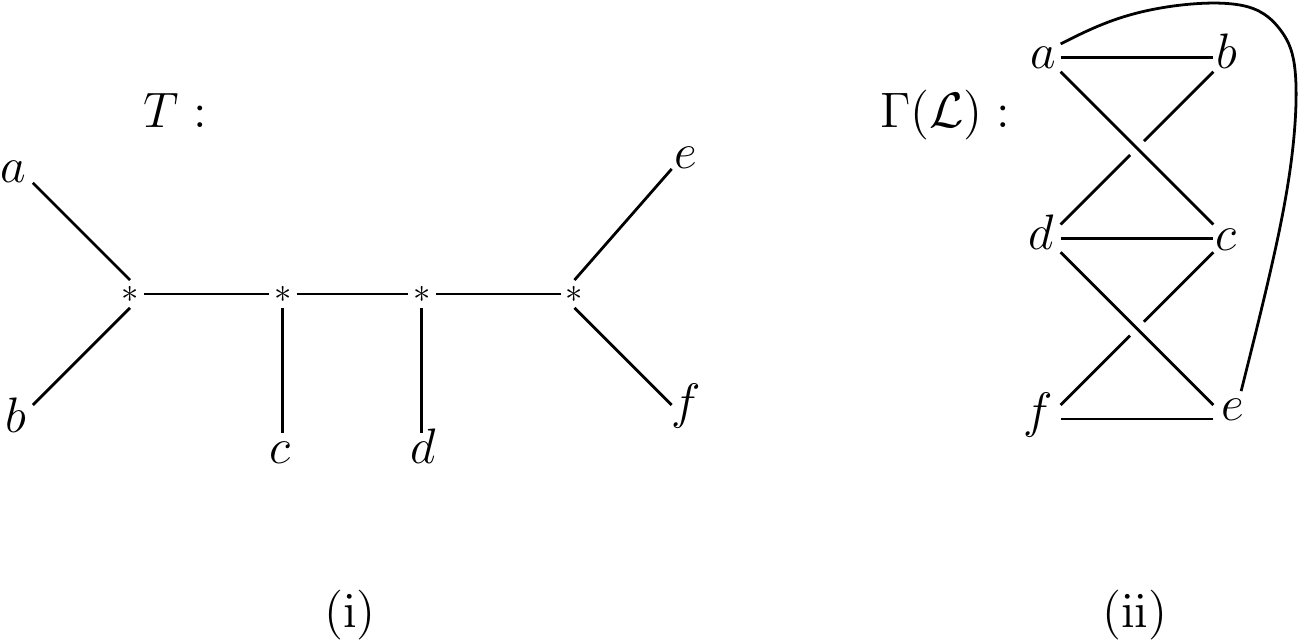}}
}
\caption{\label{figure3}}{(i) An $X-$tree $T=(V,E)$ for $X:=\{a,b,c,d,e,f\}$ with only proper $T-$cherries, and (ii) a bipartite topological lasso $\cl$ 
of $T-$rank $8 (=|E|-1)$. 
Note that both $T-$cherries $ab$ and $ef$ have a non-empty intersection with both parts 
$\{a,d,f\}$ and $\{b,c,e\}$ of the bipartition of $X$ induced by $\cl$.}
\end{center}
\end{figure}

\section{Two concluding comments}

\subsection{$X-$trees $T$ for which $\rM(T)$ is a non-binary matroid}

Let us note finally that the matroid $\rM(T)$ associated to
an $X-$tree $T$ is a non-binary matroid whenever there exist three disjoint $2$-subsets in $X$ each of which forms a $T-$cherry:
Indeed, assume that $x_1,x_2,x_3,x_4,x_5,x_6$ are six distinct elements in $X$ 
such that each of the three  pairs $x_i,x_{i+3} \,\,(i=1,2,3)$ forms a 
$T-$cherry. It is then easy to check that the two subsets $\cl_1:=\{x_1x_2, x_2x_3,x_3x_4,x_4x_5,x_5x_6,x_6x_1\}$ and $\cl_2:=\{x_1x_3, x_3x_4,x_4x_6,x_6x_1\}$ of $\ch$ are circuits in $\rM(T)$ while their symmetric difference 
$$
\cl:=\cl_1\triangle \cl_2=
\big\{x_1x_2,x_2x_3, x_3x_1,x_4x_5,x_5x_6,x_6x_4\big\}
$$ 
is an independent subset of $\ch$ in $\rM(T)$.

Clearly, this implies that a binary $X-$tree $T$ for which $\rM(T)$ is a binary matroid must be a caterpillar tree.   
More generally, an arbitrary $X-$tree $T$ for which $\rM(T)$ is a binary matroid must be either a star tree with at most five leaves or
an $X-$tree for which -- as in the case of the (binary) caterpillar trees -- two interior vertices $u$ and $v$ of $T$ exist for which the path from $u$ to $v$ in $T$ passes every interior vertex of $T$ and all of these except perhaps $u$ and $v$ have degree $3$ while the two vertices $u$ and $v$ have degree $3$ or $4$.  We will show in a separate paper that, conversely, $M(T)$ is a binary matroid whenever this holds.

\subsection{Minimal strong lassos do not form a matroid}

Although the minimal edge-weight lassos for any $X-$tree form a matroid defined on $\ch$, the same is not always true for the 
minimal strong lassos. 

To see this, let $X=\{a,b,c,d,e,f\}$ and consider the sets: 
$$\cl_1:=
\big\{ ab, ac, ad, bc, bd, cd, ef, ae, be, ce, de, df \big\},$$
and
$$\cl_2:=\big\{ab, ac, ad, bc, bd, cd, ef, ae, be, cf, df \big\}.$$

Both $\cl_1$ and $\cl_2$ are minimal strong lassos for the $X-$tree that has one interior vertex adjacent to $a,b,c,d$, and a second interior vertex adjacent to $e,f$; however,
$\cl_1$ has one more element than $\cl_2$.

\section*{{\sc Acknowledgements}}
We thank Geoffrey Whittle for helpful comments and references to some relevant matroid-theory literature.
A.D.  thanks the CAS and the MPG for financial support;   M.S. thanks the Royal Society of New Zealand under its Marsden and James Cook Fellowship scheme.

\end{document}